\documentclass[a4paper,11pt, reqno]{amsart}
\usepackage[utf8x]{inputenc}
\usepackage{amsmath, amssymb, amsthm, amsfonts}
\usepackage{fullpage}
\usepackage{times}
\usepackage{times,mathptmx}
\usepackage{tikz}
\usepackage{mathtools}
\usepackage{fontenc}
\usepackage{enumerate}
\usepackage{yfonts}

\title{Mean-Value of Product of Shifted Multiplicative Functions and Average Number of Points on Elliptic Curves.}
\author{R. Balasubramanian}
\address{Department of Mathematics\\
  Institute of Mathematical Sciences\\
  Chennai, India-600113}
\email[R. Balasubramanian]{balu@imsc.res.in}
\author{Sumit Giri}
  \address{Department of Mathematics\\
  Institute of Mathematical Sciences\\
  Chennai, India-600113}
\email[Sumit Giri]{gsumit@imsc.res.in}

\date{}
\newtheorem{theoremm}{Theorem}

\newtheorem{theorem}{Theorem}

\newtheorem{rem}{Remark}

\newtheorem{cor}{Corollary}

\newtheorem{conj}{Conjecture}

\newcommand{\N}{\mathbb{N}}

\newcommand{\Q}{\mathbb{Q}}

\newcommand{\C}{\mathbb{C}}

\newcommand{\F}{\mathbb{F}}

\begin{document}

\begin{abstract}
In this paper, we consider the mean value of the product of two real valued multiplicative functions with shifted arguments.
The functions $F$ and $G$ under consideration are close to two nicely behaved functions $A$ and $B$, such that the average value of $A(n-h)B(n)$ over any arithmetic progression is only dependent on the common difference of the progression.
We use this method on the problem of finding mean value of $K(N)$, where $K(N)/\log N$ is the expected number of primes such that a random elliptic curve over rationals has $N$ points when reduced over those primes.
\end{abstract}
 
\maketitle
\section{\textbf{Introduction}}
Let $F$ and $G$ $:$ $\N \to \C$ be non zero multiplicative functions (a function $F$ is multiplicative if $F(mn)=F(m)F(n)$ for $(m,n)=1$).
 In this paper we are interested in finding the mean value of $F(n-h)G(n)$ for a fixed integer $h$. More precisely the sum of the form 
 \begin{align}
  M_{x,h}(F,G)=\frac{1}{x}\sum_{n\le x}F(n-h)G(n).\label{eqn1}
 \end{align}
 A lot of work has been done to find the asymptotic behavior of $M_{x,h}(F,G)$ under various conditions, (see for example \cite{9}, \cite{10}, \cite{11}, \cite{12}, \cite{13}, \cite{14}).
 In many of those cases, the functions are required to be close to $1$ on the set of primes. In some cases (for example \cite{10}) convergence of suitable series involving $F$ and $G$ has been assumed.  \par
When the functions grow faster, the problem becomes more difficult. In \cite{15}, divisor function and other faster growing functions are discussed. The Euler totient function $\phi(n)$ has been studied in \cite{8} and \cite{16}. \par
In the first theorem of this paper we consider this problem for a wide class of functions with more general growth conditions. The type of functions that we consider in Theorem \ref{Thm_1} need not necessarily be multiplicative. But they can be written as
\begin{align}
F(n)=A(n)\sum_{d\mid n}f(d)\quad \text{and }\quad & G(n)=B(n)\sum_{d\mid n}g(d), \label{eqn2}\\
\intertext{where} 
 \sum_{d=1}^{\infty}\frac{|f(d)|}{d}<+\infty,\quad & \sum_{d=1}^{\infty}\frac{|g(d)|}{d}<+\infty.\label{eqn3}
\end{align}
  Further we assume the existence of two function $M(x)$ and $E_1(x)$ such that for any positive integers $a$ and $m$,
\begin{align}
 \sum_{\underset{n\equiv a(\text{mod }m)}{n\le x}}A(n-h)B(n)=\frac{1}{m}M(x)+O(E_1(x)).\label{eqn4}
\end{align}
In the first theorem we show that under the above conditions one can prove an asymptotic estimate of $M_{x,h}(F,G)$. Further in order to write the error term explicitly, we introduce two suitable monotonic functions $E_1(x)$ and $E_2(x)$ such that
\begin{align}
\sum_{d\le x}|f(d)|=O(E_2(x)), &\quad 
\sum_{d\le x}|g(d)|=O(E_3(x)).\label{eqn5}
\end{align}
 Then the first result of this paper is as follows

 \begin{theorem}\label{Thm_1}
 Let $F$ and $G$ be two arithmetic functions, satisfying (\ref{eqn2}), (\ref{eqn3}), (\ref{eqn4}) and (\ref{eqn5}) where $f$ and $g$ are multiplicative. Let $0<c<2$, such that for any large positive real number $y$, $E_i(2y)\le cE_i(y)$ for $(i=2,3)$.
 Then for any fixed positive integer $h$, 
\begin{align}
      \sum_{n\le x}F(n-h)G(n)&=C_hM(x)+O(|E_1(x)E_2(x)E_3(x)|+|\frac{M(x)}{x}(|E_2(x)|+|E_3(x)|)|),\\
     \end{align}
with $$C_h= \prod_{p}\left( 1+\sum_{j\ge 1}\frac{f(p^j)+g(p^j)}{p^j} \right )\prod_{p\mid h}\left( 1+\frac{\overset{\nu_p(h)}{\underset{i=1}{\sum}}p^iS_p(p^i)}{S_p(1)}\right)$$
where $S_p(p^i):=\underset{\min\{e_1,e_2\}=i}{\sum}\frac{f(p^{e_1})g(p^{e_2})}{p^{e_1+e_2}}$, for $i\ge 0.$
\end{theorem}

\begin{rem}
The additional condition of $f$ and $g$ being multiplicative in Theorem \ref{Thm_1} is only required to get an Euler product form of the constant $C_h$. Also note that if $\frac{F}{A}$ is multiplicative, then by m\"{o}bius inversion formula, $f$ is uniquely determined.
Also if $\frac{F}{A}$ is \textquoteleft sufficiently\textquoteright close to $1$ on primes, then (\ref{eqn3}) is satisfied for $f$. Similarly for $\frac{G}{B}$.
So for multiplicative functions the idea is to choose \textquoteleft smooth\textquoteright functions $A$ and $B$ such that $\frac{F}{A}$ and $\frac{G}{B}$ are close to $1$. 
Also $A(n-h)B(n)$ should be nicely summable on every arithmetic progression.   
\end{rem}

Before proceeding with the proof of Theorem \ref{Thm_1} we shall note down some application of the above theorem. One can directly apply it on classical Euler's totient function $\phi$ and Jordan's totient function $J_k$ . See \cite{19} and \cite{18} for more on the error term related to $\phi$ and $J_k$. Also see \cite{16} for the mean  value of the $k$-fold shifted product of $\phi$.\par

\begin{cor}\label{cor_1}
\begin{enumerate}[(a)]
 \item If $\phi(n)$ is the Euler totient function, i.e. $\phi(n)=n\underset{p\mid n}{\prod}(1-1/p)$, then for any fixed integer $h$
 \begin{align*} 
  \sum_{n\le x}\phi(n)\phi(n-h)=\frac{1}{3}x^3\prod_p(1-\frac{2}{p^2})\prod_{p\mid h}(1+\frac{1}{p(p^2-2)})+O(x^2(\log x)^2).
 \end{align*}
\item If $J_k(n)$ is the Jordan's totient function defined as $J_k(n)=n^k\underset{p\mid n}{\prod}(1-1/p^k)$, then for $k\ge 2$ and fixed integer $h$
\begin{align*}
 \sum_{n\le x}J_k(n)J_k(n-h)=\frac{x^{2k+1}}{2k+1}\prod_p(1-\frac{2}{p^{k+1}})\prod_{p\mid h}\left(1+\frac{1}{p^k(p^{k+1}-2)}\right)+O(x^{2k}).
\end{align*}
\end{enumerate}
 \end{cor}

 \par

 Proof of Corollary \ref{cor_1} follows directly from Theorem \ref{Thm_1}. In case of (a), $A(n)=B(n)=n$, while for  Jordan totient function $J_k(n)$, one takes $A(n)=B(n)=n^k$. For both the cases $f$ and $g$ can be computed using m\"{o}bius inversion.  
 
 In the next part, we discuss an application of Theorem \ref{Thm_1} in computing the mean value of the function $K(N)$ as defined in \cite{1}. Before stating the result we explain the background of this problem.

Let $E$ be an elliptic curve defined over the field of rationals $\Q$. For a primes $p$ where $E$ has good reduction, we denote by $E_p$ the reduction of $E$ modulo $p$.
Let $\F_p$ be the finite field with $p$ elements. 
Define $M_E(N)$ as 
\begin{align}
 M_E(N):=\#\{p \text{ prime}:\, E\text{ has good reduction over $p$ and } |E_p(\F_p)|=N \}. 
\end{align}
Using Hasse bound and upper bound sieve one can show that \begin{align}M_E(N)&\ll \frac{\sqrt{N}}{\log N}.\label{eqn7}\end{align} 
If $E$ has complex multiplication (CM), then Kowalski\cite{4} has shown that $$M_E(N)\ll N^\epsilon$$ for any $\epsilon >0$.\par
No stronger bound is known when $E$ is non-CM.
A naive probabilistic model suggests $M_E(N)\sim \frac{1}{\log N}$. See \cite{1} for details. Any estimate of $M_E(N)$ for a fixed $E$ is not possible. In fact using Chinese Reminder Theorem it can be shown that for giver integer $N$, the bound in (\ref{eqn7}) is attained for some $E$.
In \cite{4}, Kowalski has shown that \begin{align}\sum_{N\le x}M_E(N)=\pi(x)+O(\sqrt{x}).\label{eqn6}\end{align}

In \cite{1} David and Smith introduced an arithmetic function $K(N)$. Later they made a correction\cite{20} in the expression of $K(N)$. The corrected formula is as follows 
\begin{align}K(N):&=\prod_{p\nmid N}\left( 1-\frac{(\frac{N-1}{p})^2p+1}{(p-1)^2(p+1)}\right)\prod_{{p\mid N}}\left(1-\frac{1}{p^{\nu_p(N)}(p-1)}\right)  \end{align}
where $\nu_p$ denotes the usual $p$-adic valuation.\par 
Now let $K^*(N)=K(N)N/\phi(N)$, where $\phi(N)$ is the Euler totient function. 
In \cite{1}, David and Smith proved an asymptotic estimate for average value of $M_E(N)$ when $E$ varies over a family of curves. But their result was not unconditional. It depends on the following conjecture 
\begin{conj}[Barban–Davenport–Halberstam]\label{conj_1} Let $\theta(x; q, a) = \underset{p\le x, {p\equiv a (\text{mod }q)}}{\sum} \log p$. Let $0 < \eta \le 1$ and $\beta > 0$ be real numbers. Suppose that $X$, $Y$, and $Q$ are positive real numbers satisfying $X^\eta \le Y \le X$
and $Y /(\log X)^\beta \le Q \le Y$. Then
$$\sum_{q\le Q }\sum_{\underset{(a,q)=1}{1\le a\le q}}|\theta(X + Y ; q, a) − \theta(X; q, a) −\frac{Y}{\phi(q)}|^2\ll_{\eta,\beta} Y Q \log X.$$
 \end{conj}
\begin{rem}\label{rem_1}
For $\eta = 1$, this is the classical Barban–Davenport–Halberstam theorem. Languasco, Perelli, and Zaccagnini \cite{22} have proved the Conjecture for $\eta=\frac{7}{12}+\epsilon$, which is the best known result. Also under generalized Riemann hypothesis they could prove the conjecture for $\eta=\frac{1}{2}+\epsilon$.   
\end{rem}

Given integers $a$ and $b$, let $E_{a,b}$ be the elliptic curve defined by the Weierstrass equation
$$E_{a,b} : y^2 = x^3 + ax + b.$$

For $A$, $B > 0$, we define a set of Weierstrass equations by
$$\mathcal{C}(A, B) := \{E_{a,b} : |a|\le A, |b|\le B, ∆(E_{a,b} ) \neq 0\}.$$

In [\cite{1}, \cite{20}, \cite{21}], the following conditional result has been proved.
\begin{theoremm}\label{thmm_1}
Assume Conjecture \ref{conj_1} holds for some $\eta <\frac{1}{2}$. Let $\epsilon>0$ and $A$, $B> N^{\frac{1}{2}+\epsilon}$ such that $AB>N^{\frac{3}{2}+\epsilon}$. Then for any positive integer $R$,
\begin{align*}
\frac{1}{\#\mathcal{C}(A, B)}\sum_{E\in \mathcal{C}(A,B)}M_E(N)=\frac{K^*(N)}{\log N}+O_{\eta,\epsilon, R}(\frac{1}{(\log N)^R}). 
\end{align*}
\end{theoremm}

In order to verify the consistency of Theorem \ref{thmm_1} with unconditional results such as (\ref{eqn6}), one need to compute the mean value of $K^*(N)$ where $N\le x$ satisfies congruence conditions. For more details see  \cite{2}.\par 
In \cite{2}, Smith, Martin and Pollack have addressed this aspect. They proved that  
\begin{theoremm}\label{Thmm_2}
 For $x\ge 2$, $$\sum_{N\le x}K^*(N)=x+O(\frac{x}{(\log x)})\quad \text{and }\quad \sum_{\underset{n \text{ odd}}{N\le x}}K^*(N)=\frac{x}{3}+O(\frac{x}{(\log x)}).$$
\end{theoremm}

Using Theorem \ref{Thmm_2} and Abel's partial summation one can verify that
\begin{align*}
 \frac{1}{\#\mathcal{C}(A, B)}\sum_{E\in \mathcal{C}(A,B)}\sum_{N\le x}M_E(N)=\frac{x}{\log x}+O(\frac{x}{(\log x)^2}).
\end{align*}
So Theorem \ref{thmm_1} consistent with (\ref{eqn6}) if one consider $\pi(x)=\frac{x}{\log x}+O(\frac{x}{(\log x)^2})$.\par
But it is well known that $li(x)=\int_2^{\infty}\frac{1}{\log x} dx$ is a better approximation of $\pi(x)$ compared to $\frac{x}{\log x}$. So in order to check the consistency of Theorem \ref{thmm_1} and (\ref{eqn6}), where main term of $\pi(x)$ is taken as $li(x)$, we need significantly better bound for the error terms in Theorem \ref{Thmm_2}. In this paper we prove that. We prove 

\begin{theorem}\label{Thm_2} For $x\ge 2$,
 \begin{enumerate}[(a)] \item \begin{align*}\sum_{N\le x}K^*(N)=x+O(\log x)\end{align*} \item \begin{align*}\sum_{\underset{N \text{ odd}}{N\le x}}K^*(N)=\frac{x}{3}+O(\log x).\end{align*}\end{enumerate}
\end{theorem}

Then Theorem \ref{thmm_1} and Theorem \ref{Thm_2} together implies 
\begin{align*}
 \frac{1}{\#\mathcal{C}(A, B)}\sum_{E\in \mathcal{C}(A,B)}\sum_{N\le x}M_E(N)=li(x)+O(\frac{x}{(\log x)^R}).
\end{align*}
This provides further support to the Barban–Davenport–Halberstam conjecture. \par 
Although the function $K^*(N)$ looks like a multiplicative function it is far from it. In fact  
\begin{align}
 K^*(N)=C^*_2F^*(N-1)G^*(N)\label{eq_3}
\end{align}
where
\begin{align}
C^*_2&=\prod_{p>2}\left(1-\frac{1}{(p-1)^2}\right)\label{eq_0}\\
F^*(N)&=\prod_{\underset{p>2}{p\mid N}}\left(1-\frac{1}{(p-1)^2}\right)^{-1}\prod_{p\mid N}\left(1-\frac{1}{(p-1)^2(p+1)}\right)\label{eq_1}\\
G^*(N)&=\frac{N}{\varphi(N)}\prod_{\underset{p>2}{p\mid N}}\left(1-\frac{1}{(p-1)^2}\right)^{-1}\prod_{p\mid N}\left(1-\frac{1}{p^{\nu_p(N)}(p-1)}\right)\label{eq_2}.
\end{align}
Note that, both $F^*$ and $G^*$ are multiplicative functions. \par

In the last section of this paper we discuss the original expression of $K(N)$ as defined in [Theorem 3 ; \cite{1}]. We denote it by $\hat{K}(N)$. It was defined as follows
\begin{align}
\hat{K}(N):&=\prod_{p\nmid N}\left( 1-\frac{(\frac{N-1}{p})^2p+1}{(p-1)^2(p+1)}\right)\prod_{\underset{2\nmid \nu_p(N)}{p\mid N}}\left(1-\frac{1}{p^{\nu_p(N)}(p-1)}\right)  \prod_{\underset{2\mid \nu_p(N)}{p\mid N}}\left(1-\frac{p-\left(\frac{-N_p}{p}\right)}{p^{\nu_p(N)+1}(p-1)}\right)\label{eq_7}
\end{align}
where $\nu_p$ denotes the usual $p-$adic valuation, and $N_p:=\frac{N}{p^{\nu_p(N)}}$ denotes the $p-$free part of $N$.\par

This function cannot be written as product of two shifted multiplicative function. In \cite{2}, it is claimed that the mean of $K^*(N)$ is also equals to $1$.\par
But we show that is not true. The average turns out to be equal to $\frac{31}{30}$. Also we make improvement on the error term in the average of $\hat{K}(N)$. We prove that 

\begin{theorem}\label{Th_3}
       For $x\ge 2$, $$\sum_{N\le x}\hat{K}(N)=\frac{31}{30}x+O(\log x).$$
      \end{theorem}

The main reason behind proving this theorem separately is to show that Theorem \ref{Thm_1} can be useful in some cases where one of the shifted multiplicative functions is not multiplicative. Under suitable conditions those non-multiplicative functions can be changed to expected multiplicative form. That way Theorem \ref{Thm_1} can also be usefull in computing mean value of function.\par
In the next sections we give give proofs of the above three theorem. \par

\section{Proof of \emph{Theorem \ref{Thm_1}}}

We have  
\begin{align}
      \sum_{n\le x}F(n-h)G(n)&=\sum_{n\le x}G(n)A(n-h)\sum_{d\mid n-h}f(d)\nonumber\\
      &=\sum_{d\le x-h}f(d)\sum_{\underset{n\equiv h(\text{mod }d)}{n\le x}}G(n)A(n-h)\nonumber\\
      &=\sum_{d\le x-h}f(d)\sum_{\underset{n\equiv h(\text{mod }d)}{n\le x}}A(n-h)B(n)\sum_{d_1\mid n}g(d_1)\nonumber\\
      &=\sum_{d\le x-h}f(d)\sum_{\underset{(d,d_1)\mid h}{d_1\le x}}g(d_1)\sum_{\underset{\underset{n\equiv h(\text{mod }d)}{n\equiv 0(\text{mod }d_1)}}{n\le x}}A(n-h)B(n)\nonumber\\
      &=\sum_{d\le x-h}f(d)\sum_{\underset{(d,d_1)\mid h}{d_1\le x}}g(d_1)\left(\frac{M(x)}{[d,d_1]}+O(E_1)\right),\quad \text{where } [d,d_1]:=\text{lcm}\{d,d_1\}\nonumber\\
      &=M(x)\sum_{d\le x-h}\frac{f(d)}{d}\sum_{\underset{(d,d_1)\mid h}{d_1\le x}}\frac{g(d_1)(d,d_1)}{d_1}+O(E_1(x)E_2(x)E_3(x)).\label{eqn8}
      \end{align}

      Now, the $d$-sum and $d_1$-sum can be extended to $\infty$ to get 
      \begin{align*}
       M(x)&\sum_{d=1}^{\infty}\frac{f(d)}{d}\sum_{\underset{(d,d_1)\mid h}{d_1=1}}^{\infty}\frac{g(d_1)(d,d_1)}{d_1}\\
       \intertext{with an error term }
       O(M(x)\sum_{ 1\le d <+\infty}\frac{f(d)}{d}&\sum_{\underset{(d,d_1)\mid h}{d_1> x}}\frac{g(d_1)(d,d_1)}{d_1})  +O(M(x)\sum_{d> x-h}\frac{f(d)}{d}\sum_{\underset{(d,d_1)\mid h}{d_1\le x}}\frac{g(d_1)(d,d_1)}{d_1}).
      \end{align*}

\par Now note that
\begin{align*}
\sum_{d>x}\frac{|f(d)|}{d}&=\sum_{x<d\le 2x}\frac{|f(d)|}{d}+\sum_{2x<d\le 4x}\frac{|f(d)|}{d}+\sum_{4x<d\le 8x}\frac{|f(d)|}{d}+\cdots\\
&\le \frac{E_2(2x)}{x}+\frac{E_2(4x)}{2x}+\frac{E_2(8x)}{4x}+\cdots\\
&\le \frac{E_2(x)}{x}(c+c^2/2+c^3/4+c^4/8+\cdots)\\
&\le \frac{2c}{2-c}\frac{E_2(x)}{x}.
\end{align*}
Thus $\underset{d>x}{\sum}\frac{|f(d)|}{d}=O(\frac{E_2(x)}{x})$. Similarly $\underset{d_1>x}{\sum}\frac{|g(d_1)|}{d_1}=O(\frac{E_3(x)}{x})$.\par

Then by (\ref{eqn8})\begin{align}
                  \sum_{n\le x}F(n-h)G(n)&=M(x)\sum_{\underset{(d,d_1)\mid h}{d, d_1}}\frac{f(d)g(d_1)(d,d_1)}{dd_1}+O(|E_1(x)E_2(x)E_3(x)|+\frac{M(x)}{x}(|E_2(x)|+|E_3(x)|)). \label{eq_9*}
                    \end{align}
\par
Only thing that remains to complete the proof is to express $\underset{\underset{(d,d_1)\mid h}{d, d_1}}{\sum}\frac{f(d)g(d_1)(d,d_1)}{dd_1}$ as an Euler product.\par
To do that define the following notations
\begin{align*}
 T(p^k)&:=\frac{S_p(p^k)}{S_p(1)}=\frac{\underset{\min \{e_1,e_2\}=k}{\sum}\frac{f(p^{e_1})g(p^{e_2})}{p^{e_1+e_2}}}{\underset{\min \{e_1,e_2\}=0}{\sum}\frac{f(p^{e_1})g(p^{e_2})}{p^{e_1+e_2}}}\\
T(h_1)&:=\prod_{p\mid h_1}T(p^{\nu_p(h_1)}). 
\end{align*}
Then one can verify that $$\sum_{(d,d_1)=h_1}\frac{f(d)g(d_1)}{dd_1}=T(h_1)\sum_{(d,d_1)=1}\frac{f(d)g(d_1)}{dd_1}.$$
Now 
\begin{align*}
\underset{\underset{(d,d_1)\mid h}{d, d_1}}{\sum}\frac{f(d)g(d_1)(d,d_1)}{dd_1}&=\sum_{h_1\mid h}h_1T(h_1)\sum_{(d,d_1)=1}\frac{f(d)g(d_1)}{dd_1}\\
&=\left (\sum_{(d,d_1)=1}\frac{f(d)g(d_1)}{dd_1}\right )\prod_{p\mid h}(1+pT(p)+\cdots +p^{\nu_p(h)}T(p^{\nu_p(h)}))\\
&=\prod_{p}\left ( 1+\frac{f(p)+g(p)}{p}+\frac{f(p^2)+g(p^2)}{p^2}+\cdots \right) \prod_{p\mid h}(1+pT(p)+\cdots +p^{\nu_p(h)}T(p^{\nu_p(h)}))
\end{align*}
which proves the result.

\section{Proof of Theorem \ref{Thm_2}}

 Recall that, $$K^*(N)=C^*_2F^*(N-1)G^*(N)$$ where $C_2^*$, $F^*$ and $G^*$ are given as in (\ref{eq_0}), (\ref{eq_1}), (\ref{eq_2}).\par

  Now in this case $A(n)=B(n)=1$, hence $M(x)=x$. Also if we set  
  \begin{align}f^*(m)&=\underset{d\mid m}{\sum}\mu(d)F^*(m/d)\label{eq_17}\\ \intertext{and}
  g^*(m)&=\underset{d\mid m}{\sum}\mu(d)G^*(m/d),\label{eq_18}\end{align}
  then they are multiplicative functions. So it is enough to compute the values on prime powers. It is straight forward to check that 
 \begin{align}
 f^*(p^k)&=\left \{ \begin{array}{ll}
                         1, \quad \quad\text{if }k=0\\
                         1/(p+1)(p-2), \quad \text{if }k=1\\
                         0, \quad \quad \text{else.}\\
                        \end{array}
\right.\label{eqn9} \\
g^*(p^k)&=\left \{\begin{array}{ll}
                 1, \quad \text{if }k=0\\
                 (p-1)/p^k(p-2), \quad \text{if }k\ge 1
                \end{array}
\right.\label{eq_4}
 \end{align}
  for primes $p>2$. Also
  \begin{align}f^*(2^k)&=\left \{ \begin{array}{ll}
                         -1/3, \quad \quad\text{if }k=1\\
                         0, \quad \quad \text{if }k\ge 2\\
                        \end{array}
\right. \label{eqn10}\\
g^*(2^k)&=\left \{ \begin{array}{ll}
                 0, \quad \text{for }k= 1\\
                         1/2^{k-1}, \quad \quad \text{if }k\ge 2.\\
                        \end{array}
\right. \label{eq_5}\end{align}
         
First we shall compute the error terms. In order to do that it is enough to compute $E_1(x)$, $E_2(x)$ and $E_3(x)$ as in \emph{Theorem \ref{Thm_1}}.\par
Is is easy to see that $E_1(x)=O(1)$. \par
Now \begin{align*}
      E_2(x)&=\sum_{d\le x}|f^*(d)|\\
      &\ll \prod_{p\le x}(1+f^*(p)+f^*(p^2)+\cdots)\\
      &\ll \prod_{2<p\le x}(1+\frac{1}{(p+1)(p-2)})\\
      &=O(1).
     \end{align*}

 Also \begin{align*}
       E_3(x)&=\sum_{d_1\le x}|g^*(d_1)|\\
       &\le \prod_{p\le x}(1+g^*(p)+g^*(p^2)+\cdots)\\
       &\le \prod_{2<p\le x}(1+\frac{1}{p-2})\\
       &\ll \exp(\sum_{2<p\le x}\frac{1}{p-2})\\
       &\ll \log x.
      \end{align*}

Now only thing that remains is to compute the constant in the main term. To do that, we use the formula of $C_1$ from Theorem \ref{Thm_1}.
\par To prove (a), we use the expressions of $f^*(p^k)$ and $g^*(p^k)$ from (\ref{eqn9}), (\ref{eq_4}), (\ref{eqn10}) and (\ref{eq_5}).\par 
If $p\neq 2$
\begin{align}
 1+\sum_{i=1}^{+\infty}\frac{f^*(p^i)+g^*(p^i)}{p^i}&=1+\frac{1/(p+1)(p-2)+(p-1)/p(p-2)}{p}+\frac{p-1}{p-2}\sum_{i\ge 2}\frac{1}{p^{2i}}\nonumber \\
&=1+\frac{1}{p(p+1)(p-2)}+\frac{p-1}{p-2}\frac{1}{p^2-1}\nonumber\\
&=\frac{(p-1)^2}{p(p-2)}\nonumber\\
&=\left (1-\frac{1}{(p-1)^2}\right )^{-1}.\label{eq_21}
\end{align}

Also \begin{align}
       1+\sum_{i=1}^\infty\frac{f^*(2^i)+g^*(2^i)}{2^i}&=1+\frac{(-1/3)}{2}+\sum_{j\ge 2}\frac{1}{2^{2j-1}}=1.
\end{align}
Since $C_2^*=\underset{p>2}{\prod}\left(1-\frac{1}{(p-1)^2} \right)^{-1}$ this completes the proof of \emph{(a)}. \par
 
 To prove \emph{(b)}, we may assume that $G$ is supported on odd integers only. Hence $G(2^k)=0$ for all $k\ge 1$. In that case 
 \begin{align}
  g^*(2^k)=\left \{ \begin{array}{ll}
                 -1, \quad \text{if }k=1\\
                   0, \quad \text{if }k\ge 2.
                  \end{array}
\right.\label{eq_6}
 \end{align}
This gives
\begin{align*}
1+\sum_{i=1}^\infty\frac{f^*(2^i)+g^*(2^i)}{2^i}&=1+\frac{(-1/3)+(-1)}{2}\\
&=\frac{1}{3}.
\end{align*}
This proves \emph{(b)}.

 \section{Proof of Theorem \ref{Th_3}}

Recall that \begin{align}
              \hat{K}(N)&=C^*_2F^*(N-1)G^*_1(N),\nonumber\\
             \end{align}
where \begin{align}
       F^*(N)&=\prod_{\underset{p>2}{p\mid N}}\left(1-\frac{1}{(p-1)^2(p+1)}\right)\left( 1-\frac{1}{(p-1)^2} \right)^{-1}\nonumber\\
\intertext{and}       G^*_1(N)=\frac{N}{\phi(N)}\prod_{\underset{p>2}{p\mid N}}&\left(1-\frac{1}{(p-1)^2}\right)^{-1}\prod_{\underset{2\nmid \alpha}{p^\alpha \parallel N}}\left(1-\frac{1}{p^\alpha(p-1)} \right)\prod_{\underset{2\mid \alpha}{p^\alpha \parallel N}}\left(1-\frac{p-\left(\frac{-N_p}{p}\right)}{p^{\alpha+1}(p-1)} \right).
      \end{align}

We write $G^*_1(N)=G^*_2(N)G^*_3(N)$, where 
\begin{align}
G^*_2(N)&=\frac{N}{\phi(N)}\prod_{\underset{p>2}{p\mid N}}\left(1-\frac{1}{(p-1)^2}\right)^{-1}\prod_{\underset{2\nmid \alpha}{p^\alpha \parallel N}}\left(1-\frac{1}{p^\alpha(p-1)} \right)\nonumber\\ 
\intertext{and} 
G^*_3(N)&=\prod_{p^{2\alpha} \parallel N}\left(1-\frac{p-\left(\frac{-N_p}{p}\right)}{p^{2\alpha+1}(p-1)} \right).\label{eq_15}
\end{align}
Then $G_2^*$ is multiplicative but $G_3^*$ is not.Write 
\begin{align}
 G^*_2(N)&=\sum_{l\mid N}\hat{g}(l).\nonumber
\end{align}

Then, if $p\neq 2$, $$\hat{g}(p^k)=\left \{ \begin{array}{ll}
                                1, \text{ if }k=0\\
                                \frac{(p-1)}{p(p-2)}, \text{ if }k=1\\
                                \frac{1}{p^{2s-1}(p-2)}, \text{ if }k=2s, \, s\ge 1\\
                                -\frac{1}{p^{2s+1}(p-2)}, \text{ if }k=2s+1,\,  s\ge 1
                               \end{array}
\right.$$

and $$\hat{g}(2^k)=\left \{ \begin{array}{ll}
                                1, \text{ if }k=0\\
                                0, \text{ if }k=1\\
                                \frac{1}{2^{k-2}}, \text{ if }k=2s, s\ge 1\\
                                -\frac{1}{2^{k-1}}, \text{ if }k=2s+1, s\ge 1.
                                \end{array}
\right.$$

Our claim, which is motivated from a similar idea in \cite{2}, is that the whole computation of $\underset{N\le x}{\sum}F^*(N-1)G_1^*(N)$ remains the same even if we replace $\left(\frac{-N_p}{p}\right)$ in $G^*_3(N)$ by its expected value $0$ for every the prime other than $2$ and in case of $p=2$, we replace it by $1$. To make this rigorous, define 
\begin{align*}G^*_4(N)&=\prod_{\underset{p\neq 2}{p^{2\alpha}\parallel N}}\left(1-\frac{1}{p^{2\alpha}(p-1)}\right)\prod_{2^{2\alpha}\parallel N}\left(1-\frac{1}{2^{2\alpha+1}}\right).\end{align*}
For any $d$, $l$ with $(d,l)=1$, we claim that  
\begin{align}
 \sum_{\underset{\underset{N\equiv 0(\text{mod }l)}{N\equiv 1(\text{mod }d)}}{N\le x}}G^*_3(N)&=\sum_{\underset{\underset{N\equiv 0(\text{mod }l)}{N\equiv 1(\text{mod }d)}}{N\le x}}G^*_4(N)+O(1).\label{eqn11}
\end{align}
To prove that 

\begin{align}
 \sum_{\underset{\underset{N\equiv 0(\text{mod }l)}{N\equiv 1(\text{mod }d)}}{N\le x}}G^*_3(N)
 &=\sum_{\underset{\underset{N\equiv 0(\text{mod }l)}{N\equiv 1(\text{mod }d)}}{N\le x}}\prod_{p^{2\alpha} \parallel N}\left(1-\frac{p-\left(\frac{-N_p}{p}\right)}{p^{2\alpha+1}(p-1)} \right)\nonumber\\
 &=\sum_{\underset{\underset{N\equiv 0(\text{mod }l)}{N\equiv 1(\text{mod }d)}}{N\le x}}\prod_{p^{2\alpha} \parallel N}\left(1-\frac{1}{p^{2\alpha}(p-1)}+\frac{{\left(\frac{-N_p}{p}\right)}/p}{p^{2\alpha}(p-1)} \right).\label{eqn12}
\end{align}

From now on $l_1$, $l_2$, $l_3$ are mutually co-prime positive integers. we define the following notations
$$\psi(l_i)=\prod_{p^\beta\parallel l_i}p^\beta (p-1),$$  $$A(m, l_i)=\prod_{p\mid l_i}\frac{( \frac{-m_p}{p})}{p},$$ and $$l_3'=\prod_{p\mid l_3}p.$$
Now if $\omega(m)$ denote the number of distinct prime divisors of $m$, then with these notations, (\ref{eqn12}) is equal to 
\begin{align}
 \sum_{\underset{\underset{l_1l_2^2l_3^2\equiv 0(\text{mod }l)}{l_1l_2^2l_3^2\equiv 1(\text{mod }d)}}{l_1l_2^2l_3^2\le x}}\frac{(-1)^{\omega(l_2)}A(l_1l_2^2l_3^2, l_3)}{\psi(l_2^2l_3^2)}
& =\sum_{l_2^2l_3^2\le x}\frac{(-1)^{\omega(l_2)}}{l_3'\psi(l_2^2l_3^2)}\sum_{\underset{\underset{l_1l_2^2l_3^2\equiv 0(\text{mod }l)}{l_1l_2^2l_3^2\equiv 1(\text{mod }d)}}{l_1l_2^2l_3^2\le x}}(\frac{-l_1}{l_3'}).\label{eqn13}
\end{align}

Since $(l_1,l_3)=1$,  $(\frac{-l_1}{l'_3})$ can be replaced by $1$, for $l_3'=1, 2$, in the last summation. Also in case of other $l'_3$, the condition $(l_1,l_3)=1$ is taken care of by $(\frac{-l_1}{l'_3})$ \par
Hence (\ref{eqn13}) can be broken into two parts, namely $S(x,l,d)$ and $E_5(x)$, where

\begin{align}
S(x,l,d)&=\sum_{l_2^2\le x}\frac{(-1)^{\omega(l_2)}}{\psi(l_2^2)}\sum_{\underset{\underset{l_1l_2^2\equiv 0(\text{mod }l)}{l_1l_2^2\equiv 1(\text{mod }d)}}{l_1l_2^2\le x}}1+\sum_{\underset{(l_2,2)=1}{l_2^2 2^{2\gamma}\le x}}\frac{(-1)^{\omega(l_2)}}{2\psi(l_2^2 2^{2\gamma})}\sum_{\underset{\underset{l_1l_2^2 2^{2\gamma}\equiv 0(\text{mod }l)}{l_1l_2^2 2^{2\gamma}\equiv 1(\text{mod }d)}}{l_1l_2^2 2^{2\gamma}\le x}}1\nonumber\\
\intertext{and}
E_5(x)&=\sum_{\underset{l_3'\ge 3}{\underset{(l_2,l_3)=1}{l_2^2l_3^2\le x}}}\frac{(-1)^{\omega(l_2)}}{l_3'\psi(l_2^2l_3^2)}\sum_{\underset{\underset{l_1l_2^2l_3^2\equiv 0(\text{mod }l)}{l_1l_2^2l_3^2\equiv 1(\text{mod }d)}}{l_1l_2^2l_3^2\le x}}(\frac{-l_1}{l_3'}).\nonumber
\end{align}

If we rewrite $G^*_4$ as
$$G^*_4(n)=\prod_{\underset{p\neq 2}{p^{2\alpha}\parallel N}}\left(1-\frac{1}{p^{2\alpha}(p-1)}\right)\prod_{2^{2\alpha}\parallel N}\left(1-\frac{1}{2^{2\alpha}}+\frac{1}{2^{2\alpha+1}}\right),$$
then it is easy to check that
$$\sum_{\underset{\underset{N\equiv 0(\text{mod }l)}{N\equiv 1(\text{mod }d)}}{N\le x}}G^*_4(N)=S(x,l,d).$$
For $E_5$, note that the congruence relations in the last summation has no solution unless $(l_2l_3,d)=1$. So if solutions exists, then there exists $a_1$, $a_2\cdots a_{\phi(l_2)}$, such that the congruence conditions along with the condition $(l_1,l_2)=1$ is equivalent to any one of the following $$l_1\equiv a_i (\text{ mod }M_{d,l,l_2,l_3}),\quad i=1,2,\cdots,\phi(l_2)$$
with $(M_{d,l,l_2,l_3},l'_3)=1$.\par
Then for each fixed $a_i$, the set $\{a_i,\, a_i+M_{d,l,l_2,l_3},\, a_i+2M_{d,l,l_2,l_3}, \cdots ,\, a_i+(l'_3-1)M_{d,l,l_2,l_3}\}$ runs over all possible residue class module $l'_3$ exactly once.
Hence using the fact that $$\sum_{a=1}^{l_3'}(\frac{a}{l_3'})=0 \quad \text{for } l_3'\ge 3,$$
we get \begin{align*}
E_5(x)&=\sum_{l_2^2l_3^2\le x}\frac{(-1)^{\omega(l_2)}}{l_3'\psi(l_2^2l_3^2)}[0+O(\phi(l_2)l_3')]\nonumber\\
&=O(\sum_{\underset{(l_2,l_3)=1}{l_2^2l_3^2\le x}}\frac{l_2}{\psi(l_2^2l_3^2)})\nonumber\\
&=O(\sum_{\underset{(l_2,l_3)=1}{l_2\le \sqrt{x}}}\frac{l_2}{\psi(l_2^2)})\nonumber\\
&=O(\sum_{l_2\le \sqrt{x}}\frac{1}{\psi(l_2)})\nonumber\\
&=O(1).
\end{align*} 
Which proves the claim.

Now with these notations, where $f^*(d)$ is as in (\ref{eq_17}), we have
\begin{align}
      \sum_{N\le x}F^*(N-1)G^*_1(N)&=\sum_{N\le x}G^*_1(N)\sum_{d\mid N-1}f^*(d)\nonumber\\
      &=\sum_{d\le x-1}f^*(d)\sum_{\underset{N\equiv 1(\text{mod }d)}{N\le x}}G^*_1(N)\nonumber\\
      &=\sum_{d\le x-1}f^*(d)\sum_{\underset{N\equiv 1(\text{mod }d)}{N\le x}}G^*_2(N)G^*_3(N)\nonumber\\
      &=\sum_{d\le x-1}f^*(d)\sum_{\underset{N\equiv 1(\text{mod }d)}{N\le x}}G^*_3(N)\sum_{l\mid N}\hat{g}(l)\nonumber\\
      &=\sum_{d\le x-1}f^*(d)\sum_{\underset{(l,d)=1}{l\le x}}\hat{g}(l)\sum_{\underset{\underset{N\equiv 0(\text{mod }l)}{N\equiv 1(\text{mod }d)}}{N\le x}}G^*_3(N).\nonumber\\
      \intertext{Now, using (\ref{eqn11}) we get}
   \sum_{N\le x}F^*(N-1)G^*_1(N)   &=\sum_{d\le x-1}f^*(d)\sum_{\underset{(l,d)=1}{l\le x}}\hat{g}(l)\sum_{\underset{\underset{N\equiv 0(\text{mod }l)}{N\equiv 1(\text{mod }d)}}{N\le x}}G^*_4(N)\nonumber+O(\sum_{d\le x-1}|f^*(d)|\sum_{\underset{(l,d)=1}{l\le x}}|\hat{g}(l)|)\\
      &=\sum_{d\le x-1}f^*(d)\sum_{\underset{N\equiv 1(\text{mod }d)}{N\le x}}G^*_2(N)G^*_4(N)+O(\log x)\nonumber\\
      &=\sum_{N\le x}F^*(N-1)G^*_2(N)G^*_4(N)+O(\log x).
      \end{align}

To compute the main term, note that if $G^*_2(N)G^*_4(N)=\underset{l\mid n}{\sum}g_1^*(l)$, then

\begin{align*}
g^*_1(p^k)&=\left \{\begin{array}{ll}
                 1, \quad \text{if }k=0\\
                 (p-1)/p^k(p-2), \quad \text{if }k\ge 1
                \end{array}
\right.\\ 
\intertext{and}
g^*_1(2^k)&=\left \{ \begin{array}{ll}
               0, \text{ if }k=2s-1, s\ge 1\\
               \frac{3}{2^{2s}}\text{ if }k=2s, s\ge 1.
              \end{array}
\right.\end{align*}

So in order to compute the constant in the main term it is enough to compute $(1+\frac{f^*(2)+g_1^*(2)}{2}+\frac{f_1^*(2^2)+g_1^*(2^2)}{2^2}+\cdots)$, because other factors corresponding to the primes $p(\neq 2)$ cancels out with the constant $C^*_2=\underset{p>2}{\prod}\left(1-\frac{1}{(p-1)^2} \right)$.\\
Now
\begin{align*}
(1+\frac{f^*(2)+g_1^*(2)}{2}+\frac{f^*(2^2)+g_1^*(2^2)}{2^2}+\cdots)&=(1+\frac{(-1/3)}{2}+\frac{3/2^2}{2^2}+\frac{3/2^4}{2^4}+\frac{3/2^6}{2^6}+\cdots)\\ 
&=(1-\frac{1}{6}+\frac{1}{5})\\
&=\frac{31}{30}.
\end{align*}

 \begin{rem}
  In \cite{16}, Mirsky gave a proof of part (a) of Corollary \ref{cor_1}. In that same paper he discussed how to approach the problem of $k$-fold product.
  More precisely summation of the form $\sum_{n\le x}f_1(n-h_1)f_2(n-h_2)\cdots f_k(n-h_k)$, where each of $f_i$ are multiplicative.
  \end{rem}

\end{document}